
\def\smallprint{\magnification=\magstep0
     \hsize6.5truein
     \vsize9truein
     \count17=2}

\def\discrversionA#1#2{\ifnum\count17=1#1\else#2\fi}
\def\discrversionB#1#2{\ifnum\count17=2#1\else#2\fi}
\def\discrversionC#1#2{\ifnum\count17=3 #1
                   \else #2\fi}
\def\discrversionD#1#2{\ifnum\count17=4#1\else#2\fi}
\def\discrversionE#1#2{\ifnum\count17=5#1\else#2\fi}
\def\filename#1{\ifnum\count17=1 \noindent filename #1\else\fi}
\def\versiondate#1{\ifnum\count17=1

               \medskip

               \hfill Version of #1

               \medskip
               \else\ifnum\count17=2

               \medskip

               \hfill Version of #1

               \medskip
               \else  \fi\fi}
\def\today{\number\day.\number\month.\number\year}
\def\nversiondate{\ifnum\count17=1

              \medskip

              \hfill Version of \today
                \else\ifnum\count17=2

              \medskip

              \hfill Version of \today

              \medskip
                \else \fi\fi}
\def\discrpage{\ifnum\count17=3{\vfill\eject}
           \else\ifnum\count17=5{\vfill\eject}
           \else\bigskip\fi\fi}

\def\versionref{\ifnum\count17=1\workingref
           \else\ifnum\count17=2\workingref
           \else\printerref\fi\fi}
\def\versionhead{\ifnum\count17=3 \hfill PRINTER'S VERSION

                        \bigskip
                   \else\ifnum\count17=4 \hfill REFEREE'S VERSION

                        \bigskip
                   \else \fi\fi}

\def\restr{\hbox{$\hskip-0.05em\restriction\hskip-0.05em$}}

\def\frakc{\frak c\phantom{.}}

\def\loibr{\enskip]}

\documentstyle{amsppt}


\def\workingref{
\def\BSMHF{[BSM86]}
\def\DiHD{[Di84]}
\def\DUGG{[DU77]}
\def\DrHF{[Dr86]}
\def\GoIJ{[Go90]}
\def\ITFG{[IT67]}
\def\LNFG{[LN67]}
\def\LTGG{[LT77]}
\def\McSHC{[McS83]}
\def\RuGD{[Ru74]}
\def\SFpIJ{[SFp90]}
\def\TaHD{[Ta84]}
\def\TaHG{[Ta87]}
}

\def\printerref{
\def\BSMHF{[1]}
\def\DiHD{[2]}
\def\DUGG{[3]}
\def\DrHF{[4]}
\def\GoIJ{[5]}
\def\ITFG{[6]}
\def\LNFG{[7]}
\def\LTGG{[8]}
\def\McSHC{[9]}
\def\RuGD{[10]}
\def\SFpIJ{[11]}
\def\TaHD{[12]}
\def\TaHG{[13]}
}

\smallprint

\versionref
\versionhead
\filename{FMp91.tex}

\filename{FMp91a.tex}

\versiondate{20.1.92}

\medskip

\centerline{\bf On the integration of
vector-valued functions}

\centerline{\smc D.H.Fremlin \& J.Mendoza}

\centerline{\it University of Essex, Colchester, England}

\centerline{\it Universidad Complutense de Madrid, Madrid, Spain}

\medskip

{\narrower\narrower\smallskip
\noindent We discuss
relationships between the McShane, Pettis, Talagrand and Bochner
integrals.
\smallskip}

\medskip

\noindent{\bf Introduction} A large number of different methods of
integration
of
Banach-space-valued functions have been introduced, based on the various
possible constructions of the Lebesgue integral.   They commonly run
fairly closely together when the range space is separable (or has
w$^*$-separable dual) and diverge more or less sharply for general range
spaces.   The McShane integral as described by \GoIJ\ is derived from
the \lq gauge-limit' integral of \McSHC.   Here we
give both positive and negative results concerning it
and the other
three integrals listed above.

\medskip

{\bf 1A Definitions} We recall the following definitions.   Let
$(S,\Sigma,\mu)$
be a probability space and $X$ a Banach space, with dual $X^*$.

\medskip

{\bf (a)} A function $\phi:S\to X$ is {\bf Pettis integrable} if for
every $E\in\Sigma$ there is a $w_E\in X$ such that $\int_E
f(\phi(x))\mu(dx)$ exists and is equal to $f(w_E)$ for every $f\in X^*$;
in this case $w_S$ is the {\bf Pettis integral} of $\phi$, and the map
$E\mapsto w_E:\Sigma\to X$ is the {\bf indefinite Pettis integral} of
$\phi$.

\medskip

{\bf (b)} A function $\phi:S\to X$ is {\bf Talagrand integrable}, with
{\bf Talagrand integral} $w$, \linebreak if $w = \lim_{n\to\infty}{1\over
n}\sum_{i<n}\phi(s_i)$ for almost all sequences $\langle
s_i\rangle_{i\in\Bbb N}\in S^{\Bbb N}$, where $S^{\Bbb N}$ is given its
product probability.   (See \TaHG, Theorem 8.)

\medskip

{\bf (c)} A function $\phi:S\to X$ is {\bf Bochner integrable}, with
{\bf Bochner integral} $w$, if for every $\epsilon>0$ we can find a
partition $E_0,\ldots,E_n$ of $S$ into measurable sets and vectors
$x_0,\ldots,x_n\in X$ and an integrable function $h:S\to\Bbb R$ such
that
$\int h\le\epsilon$, $\|\phi(t)-x_i\|\le h(t)$ for $t\in E_i$, $i\le n$ and
$\|w-\sum_{i\le n}\mu E_i.x_i\|\le\epsilon$.

\medskip

{\bf 1B} Now we come to an integral which has been
defined for functions with domains which are intervals in $\Bbb R$.
In fact it can be satisfactorily generalized to very much wider
contexts;  but as the extension involves ideas from topological measure
theory unnecessary for the chief results of this paper, we confine
ourselves here to the original special case.

\medskip

\def\vtmp{

\centerline{$\langle([a_i,b_i],t_i)\rangle_{i\le n}$ 
}

\noindent}

\noindent{\bf Definitions}
A {\bf McShane partition} of $[0,1]$ is a finite sequence
\discrversionC{\vtmp}{$\langle([a_i,b_i],t_i)\rangle_{i\le n}$ 
} such that $\langle
[a_i,b_i]\rangle_{i\le n}$ is a
non-overlapping family of  intervals
covering $[0,1]$ and $t_i\in
[0,1]$ for each $i$.   A {\bf gauge} on $[0,1]$ is a function
$\delta:[0,1]\to\loibr 0,\infty[$.
A McShane partition
$\langle([a_i,b_i],t_i)\rangle_{i\le n}$ is {\bf subordinate} to a gauge
$\delta$ if $t_i-\delta(t_i)\le a_i\le b_i\le t_i+\delta(t_i)$ for every
$i\le n$.

Now let $X$ be a Banach space.   Following \GoIJ, we say that a function
$\phi:[0,1]\to
X$ is {\bf McShane integrable}, with {\bf
McShane integral} $w$, if for every $\epsilon>0$ there is a gauge
$\delta:[0,1]\to\loibr 0,\infty[$ such that
$\|w-\sum_{i\le n}(b_i-a_i)\phi(t_i)\|\le\epsilon$ for every McShane
partition $\langle([a_i,b_i],t_i)\rangle_{i\le n}$ of $[0,1]$
subordinate to $\delta$.

\medskip

{\bf 1C Summary of results} With four integrals to play with, a good
many questions can be asked;  and the situation is complicated by the
fact that certain natural restrictions which may be put on the space $X$
and the function $\phi$ change the answers.   We therefore set out the
facts in a
semi-tabular form.   We give references to the literature for those
which are already known, and references to paragraphs below to those
which we believe to be new.

\medskip

{\bf (a)}   Consider first the situation in which no restriction is
placed on the Banach space $X$ nor on the function $\phi:[0,1]\to X$.
In this context, a Bochner integrable function is Talagrand
integrable, a Talagrand integrable function is Pettis integrable, and
the integrals coincide whenever defined (\TaHG, Theorem 8).
A Bochner integrable function is McShane integrable (\GoIJ, Theorem 16)
(in fact, a measurable Pettis integrable function is McShane integrable
-- see \GoIJ, Theorem 17);   a McShane integrable function is Pettis
integrable (2C below).

None of the implications here can be reversed.   To see this, it is
enough to find a McShane integrable function which is not Talagrand
integrable (3A, 3G) and a Talagrand integrable function which
is not McShane integrable (3E).

\medskip

{\bf (b)}  Now suppose that the unit ball $B_1(X^*)$ of the dual $X^*$
of $X$ is separable and that $\phi$ is bounded.   In this case, a
McShane integrable function must be Talagrand integrable (2M).   We
ought to observe at this point that in these circumstances the continuum
hypothesis, for instance, is enough to make any {\it Pettis} integrable
function Talagrand integrable (\TaHD, 6-1-3), and that it remains
conceivable that this is a theorem of ZFC (see \SFpIJ).   But our result
in 2M does not depend on any special axiom.

In this context it is still true that a McShane integrable function need
not be Bochner integrable (3F) and that a Talagrand integrable
function need not be McShane integrable (3E).

\medskip

{\bf (c)}  If we take $X$ to be separable, but allow $\phi$ to be
unbounded, then the Bochner and Talagrand integrals coincide (see
2K below),
and the McShane and Pettis integrals coincide (2D).

There is still a McShane integrable function which is not Talagrand
integrable (3G).

\medskip

{\bf (d)}  For separable $X$ and bounded $\phi$, the Bochner and Pettis
integrals coincide (2K), so all four integrals here coincide.

\medskip

{\bf (e)}  Finally, the same is true, for unbounded $\phi$, if $X$ is
finite-dimensional.

\discrpage

\filename{FMp91b.tex}

\versiondate{20.1.92}

\noindent{\bf 2. Positive results}
In this section we give our principal positive results.  A McShane
integrable function is Pettis integrable (2C);  using this we are able
to prove a convergence theorem for McShane integrable functions (2I)
with some corollaries (2J).   We conclude by showing that a bounded
McShane
integrable function from $[0,1]$ to a space with
w$^*$-separable dual unit ball is Talagrand integrable (2M).

\medskip

{\bf 2A} We approach the first result by means of the \lq Dunford
integral'.   Recall that a function $\phi:[0,1]\to X$ is {\bf Dunford
integrable} if $h\phi:[0,1]\to\Bbb R$ is integrable for every $h\in
X^*$;
in this case we have an {\bf indefinite Dunford integral} $\nu:\Sigma\to
X^{**}$, where $\Sigma$ is the algebra of Lebesgue measurable subsets of
$[0,1]$, given by the formula $(\nu E)(h)=\int_E h\phi$ for every $h\in
X^*$, $E\in\Sigma$ (\TaHD, 4-4-1 or \DUGG, p.~52, Lemma 1).   
Thus a Pettis integrable function is just
a Dunford integrable function whose indefinite integral takes values in
$X$ (identified, of course, with its canonical image in $X^{**}$).   Now
we have the following general characterization of Pettis integrable
functions on $[0,1]$ (which in part, at least, is already known;  cf.
\DrHF, 4.2).

\medskip

{\bf 2B Proposition} Let $X$ be a Banach space and $\phi:[0,1]\to X$ a
Dunford integrable function with indefinite integral $\nu:\Sigma\to
X^{**}$.   Suppose that $\nu([a,b])\in X$ for every subinterval $[a,b]$
of $[0,1]$.   Then the following are equivalent:

\quad(i) $\phi$ is Pettis integrable;

\quad(ii) for every sequence $\langle[a_i,b_i]\rangle_{i\in\Bbb N}$ of
non-overlapping subintervals of $[0,1]$,  $\sum_{i\in\Bbb
N}\nu([a_i,b_i])$ exists in $X$ (for the norm of $X$);

\quad(iii) for every $\epsilon>0$ there is an $\eta>0$ such that $\|\nu
E\|\le\epsilon$ whenever $\mu E\le\eta$;

\quad(iv) $\nu$ is countably additive.

\medskip

\noindent{\bf proof (i)$\Rightarrow$(ii)} is a theorem of Pettis (see
\DUGG, II.3.5).

\medskip

{\bf (ii)$\Rightarrow$(iii)}  The point is that $\{h\phi:h\in
B_1(X^*)\}$ is uniformly integrable.   To see this, it is enough to show
that

\centerline{$\lim_{n\to\infty}\sup\{\int_{G_n}|h\phi|:h\in
B_1(X^*)\}=0$}

\noindent for every disjoint sequence $\langle G_n\rangle_{n\in\Bbb N}$
of open sets in $[0,1]$ (see \DiHD, VII.14).   But given such a
sequence, and a sequence $\langle h_n\rangle_{n\in\Bbb N}$ in
$B_1(X^*)$, set $\alpha_n=\int_{G_n}|h_n\phi|$ for each $n$.   Then we
can find for each $n$ a set $F_n\subseteq G_n$, a finite union of closed
intervals, such that $\int_{F_n}|h_n\phi|\ge
\alpha_n-2^{-n}$.   Now we can find a sequence $\langle
[a_i,b_i]\rangle_{i\in\Bbb N}$ of
non-overlapping intervals, and an increasing sequence $\langle
k_n\rangle_{n\in\Bbb N}$ of integers, such that $F_n=\bigcup_{k_n\le
i<k_{n+1}}[a_i,b_i]$ for each $n$.   So

\centerline{$\alpha_n-2^{-n}\le\|\sum_{k_n\le
i<k_{n+1}}\nu([a_i,b_i])\|$}

\noindent for each $n$.   But as $\sum_{i\in\Bbb N}\nu([a_i,b_i])$
exists in $X^{**}$, $\lim_{n\to\infty}\alpha_n=0$, as required.

Now it follows that for every $\epsilon>0$ there is an $\eta>0$ such
that $\int_E|h\phi|\le\epsilon$ whenever $h\in B_1(X^*)$, $E\in\Sigma$
and $\mu E\le\eta$, writing $\mu$ for Lebesgue measure;  so that
$\|\nu E\|\le\epsilon$ whenever $\mu E\le\eta$.

\medskip

{\bf (iii)$\Rightarrow$(iv)} is elementary.

\medskip

{\bf (iv)$\Rightarrow$(i)}
Our original hypothesis was that $\nu E\in X$ for intervals $E$;  it
follows that
$\nu E\in X$ whenever $E$ is a finite union of intervals.   Because
$\nu$ is countably additive, $\nu E\in X$ whenever $E$ is open, and
therefore whenever $E$ is G$_{\delta}$;  but also of course $\nu E=0\in
X$ whenever $\mu E=0$, so $\nu E\in X$ for every $E\in\Sigma$, and
$\phi$ is Pettis integrable.

\medskip

{\bf 2C Theorem} Let $X$ be a Banach space and $\phi:[0,1]\to X$ a
McShane integrable function.   Then $\phi$ is Pettis integrable.

\medskip

\noindent{\bf proof} As remarked in \GoIJ, Theorem 8, $\phi$ is Dunford
integrable;  let $\nu:\Sigma\to X^{**}$ be its indefinite Dunford
integral.   We know also that $\nu([a,b])\in X$ for every subinterval
$[a,b]$ of $[0,1]$ (\GoIJ, Theorem 4).   So we seek to show that (ii) of
2B above holds true.

Let $\epsilon>0$.
Let $\delta:[0,1]\to\loibr 0,\infty[$ be a gauge such that $\|\nu([0,1])
-\sum_{i\le n}(b_i-a_i)\phi(t_i)\|\le\epsilon$ whenever
$\langle([a_i,b_i],t_i)\rangle_{i\le n}$ is a McShane partition of
$[0,1]$ subordinate to $\delta$.   Fix a particular McShane partition
$\langle([a_i,b_i],t_i)\rangle_{i\le n}$ of $[0,1]$ subordinate
to $\delta$, and set $M=\sup_{i\le n}\|\phi(t_i)\|$.
We claim that if $E\subseteq[0,1]$ is a finite union of closed intervals
then $\|\nu E\|\le M\mu E+2\epsilon$.   To see this, express $E$ as
$\bigcup_{j\le m}[c_j,d_j]$ where the $[c_j,d_j]$ are
non-overlapping, and let $\eta>0$.   For each $i\le n$, we can express
$[a_i,b_i]\setminus\text{int }E$ as a (possibly empty) finite union
of
non-overlapping intervals $[a_{ik},b_{ik}]$ for $k<r(i)$;  write
$t_{ik}=t_i$ for $i\le n$, $k<r(i)$.   Then we see that

\centerline{$\|\sum_{i\le n}(b_i-a_i)\phi(t_i)
-\sum_{i\le n}\sum_{k<r(i)}(b_{ik}-a_{ik})\phi(t_{ik})\|
\le M\mu E$.}

\noindent Next, for each $j\le m$ we can find a McShane partition
$\langle([c_{jk},d_{jk}],u_{jk})\rangle_{k\le q(j)}$ of $[c_j,d_j]$,
subordinate to $\delta$, such that

\centerline{$\|\nu([c_j,d_j])-\sum_{k\le q(j)}(d_{jk}-
c_{jk})\phi(u_{jk})\|\le\eta$.}

\def\vtmp{

\centerline{$\langle([a_{ik},b_{ik}],t_{ik})\rangle_{i\le n,k<r(i)}$,
$\langle([c_{jk},d_{jk}],
u_{jk})\rangle_{j\le m,k\le q(j)}$ 
}

\noindent together}

\noindent Assembling these, we see that
\discrversionC{\vtmp}{
$\langle([a_{ik},b_{ik}],t_{ik})\rangle_{i\le n,k<r(i)}$ and
$\langle([c_{jk},d_{jk}],
u_{jk})\rangle_{j\le m,k\le q(j)}$ 
together}
form a McShane partition of $[0,1]$ subordinate to $\delta$, so that

\centerline{$\|\nu([0,1])
-\sum_{i\le n,k<r(i)}(b_{ik}-a_{ik})\phi(t_{ik})
-\sum_{j\le m,k\le q(j)}(c_{jk}-d_{jk})\phi(t_{jk})\|\le\epsilon$.}

\noindent
Also, of course,

\centerline{$\|\nu([0,1])
-\sum_{i\le n}(b_i-a_i)\phi(t_i)\|\le\epsilon$.}

\noindent Putting these formulae together, we get

\centerline{$\|\sum_{j\le m}\nu([c_j,d_j])\|\le M\mu
E+2\epsilon+m\eta$.}

\noindent But $\eta$ was arbitrary, so $\|\nu E\|\le M\mu E+2\epsilon$,
as claimed.

Now this means that if $\langle[c_j,d_j]\rangle_{j\in\Bbb N}$ is any
sequence of
non-overlapping intervals,

\centerline{$\limsup_{m\to\infty}\sup_{l\ge m}\|\sum_{m\le j\le
l}\nu([c_j,d_j])\|\le 2\epsilon$.}

\noindent But of course $\epsilon$ was arbitrary, so the limit must
actually be $0$;  accordingly $\sum_{j\in\Bbb N}\nu([c_j,d_j])$ is
defined.

Thus (ii) of 2B is satisfied, and $\phi$ is Pettis integrable.

\medskip

{\bf 2D Corollary} Let $X$ be a separable Banach space.   Then a
function $\phi:[0,1]\to X$ is McShane integrable iff it is Pettis
integrable.

\medskip

\noindent{\bf proof} If $\phi$ is Pettis integrable, it is measurable,
because $X$ is separable;  so it is McShane integrable by Theorem 17 of
\GoIJ.   Now 2C gives the reverse implication.

\medskip

{\bf 2E} As a further consequence of 2C we have the following.

\medskip

\noindent{\bf Theorem}  Let $X$ be a Banach space and $\phi:[0,1]\to
X$ a McShane integrable function.   Then for any measurable $E\subseteq
[0,1]$ the function $\phi_E=\phi\times\chi(E):[0,1]\to X$, defined by
writing
$\phi_E(t)=\phi(t)$ if $t\in E$ and $0$ otherwise, is McShane
integrable.

\medskip

\noindent{\bf proof} Again write $\nu$ for the indefinite integral of
$\phi$;  we now know that $\nu$ takes all its values in $X$.   Let
$\epsilon>0$.   By 2B, we can find an $\eta>0$ such that $\|\nu
H\|\le\epsilon$ whenever $\mu H\le\eta$.   Let $F_1\subseteq E$,
$F_2\subseteq [0,1]\setminus E$ be closed sets such that $\mu F_1+\mu
F_2\ge 1-\eta$.   Let $\delta:[0,1]\to\loibr 0,\infty[$ be a gauge such
that $\|\nu([0,1])-\sum_{i\le n}(b_i-a_i)\phi(t_i)\|\le\epsilon$
whenever $\langle([a_i,b_i],t_i)\rangle_{i\le n}$ is a McShane partition
subordinate to $\delta$.   Let $\delta_1:[0,1]\to\loibr 0,\infty[$ be
such that $\delta_1(t)\le\delta(t)$ for every $t$ and
$[t-\delta_1(t),t+\delta_1(t)]\cap F_j=\emptyset$ whenever $t\notin
F_j$, for both $j\in\{1,2\}$.

 Suppose that $\langle([a_i,b_i],t_i)\rangle_{i\le n}$ is a McShane
partition of $[0,1]$ subordinate to $\delta_1$.   We seek to estimate
$\|\nu E-\sum_{i\le n}(b_i-a_i)\phi_E(t_i)\|$.
Set $I=\{i:i\le n,\,t_i\in E\}$, $H=\bigcup_{i\in I}[a_i,b_i]$,
Then we must have $F_1\subseteq H\subseteq[0,1]\setminus F_2$.
Consequently $\mu(H\triangle E)\le\eta$ and
$\|\nu E-\nu H\|\le 2\epsilon$.   But now recall that by Theorem 5 of
\GoIJ\ we know that
$\|\nu H-\sum_{i\in I}(b_i-a_i)\phi(t_i)\|\le\epsilon$.
So

\centerline{$\|\nu E-\sum_{i\le n}(b_i-a_i)\phi_E(t_i)\|
=\|\nu E-\sum_{i\in I}(b_i-a_i)\phi(t_i)\|\le
2\epsilon$.}

\noindent As $\epsilon$ is arbitrary, $\phi_E$ is McShane integrable.

\medskip

{\bf 2F} For the next theorem of this section, we need to recall
some
well-known facts concerning vector measures.   Suppose that $\Sigma$ is
a
$\sigma$-algebra of sets and $X$ a Banach space.

\medskip

{\bf (a)} Let us say that a function $\nu:\Sigma\to X$ is \lq weakly
countably
additive' if $f(\nu(\bigcup_{i\in\Bbb N}E_i))=\sum_{i\in\Bbb N}f(\nu
E_i)$ for every disjoint sequence $\langle E_i\rangle_{i\in\Bbb N}$ in
$\Sigma$ and every $f\in X^*$.   The first fact is that in this case
$\nu$ is countably additive, that is, $\sum_{i\in\Bbb N}\nu E_i$ is
unconditionally summable to $\nu(\bigcup_{i\in\Bbb N}E_i)$ for the norm
topology whenever $\langle E_i\rangle_{i\in\Bbb N}$ is a disjoint
sequence of measurable sets with union $E$ (`Orlicz-Pettis theorem',
\TaHD, 2-6-1 or \DUGG, p.~22, Cor. 4).

\medskip

{\bf (b)} If now
$\mu$ is a measure with domain $\Sigma$ such that $\nu E=0$ whenever
$\mu E=0$, then for every $\epsilon>0$ there is a $\delta>0$ such that
$\|\nu E\|\le\epsilon$ whenever $\mu E\le\delta$.
(\DUGG, p.~10, Theorem 1.)

\medskip

{\bf (c)} Thirdly,
suppose
that
$\langle\nu_n\rangle_{n\in\Bbb N}$ is a sequence of countably
additive functions from $\Sigma$ to $X$ such that
$\nu E=\lim_{n\to\infty}\nu_n E$ exists in $X$, for the weak topology of
$X$, for every $E\in\Sigma$;
then $\nu$ is countably
additive.   (Use Nikod\'ym's theorem (\DiHD, p.~90) to see that $\nu$
is weakly countably additive.)

\medskip

{\bf 2G Lemma} Let
$X$ be a Banach space.   If $\phi:[0,1]\to X$ is
McShane integrable with McShane integral $w$, then

\centerline{$\|w\|\le\underline{\int}\|\phi(t)\|\mu(dt)$.}

\medskip

\noindent{\bf proof} Take any $f$ in the unit ball of $X^*$.   By \GoIJ,
Theorem 8,
$f(w)$ is the McShane integral of $f\phi:[0,1]\to\Bbb R$,   and by
6-4 and 6-5 of \McSHC\
this is the ordinary integral of $f\phi$.   So we have

\centerline{$|f(w)|=|\int
f\phi|\le\int|f\phi|\le\underline{\int}\|\phi\|$.}

\noindent As $f$ is arbitrary, $\|w\|\le\underline{\int}\|\phi\|$.

\medskip

{\bf 2H Lemma} Let $X$ be a Banach space and $\phi:[0,1]\to X$ a
McShane integrable function;  let $\epsilon>0$.   Then there is a gauge
$\delta:[0,1]\to\loibr 0,\infty[$ such that $\|\int_E\phi
-\sum_{i\le n}\mu E_i\phi(t_i)\|\le\epsilon$ whenever $E_0,\ldots,E_n$
are disjoint measurable subsets of $[0,1]$, $t_0,\ldots,t_n\in[0,1]$ and
$E_i\subseteq[t_i-\delta(t_i),t_i+\delta(t_i)]$ for every $i$.

\medskip

\noindent{\bf proof} Let $\delta$ be a gauge such that $\|\int\phi
-\sum_{i\le n}(b_i-a_i)\phi(t_i)\|\le\epsilon$ whenever
$\langle([a_i,b_i],t_i)\rangle_{i\le n}$ is a McShane partition of
$[0,1]$ subordinate to $\delta$.   Let $E_0,\ldots,t_n$ be as in the
statement of the lemma;  set $M=\max_{i\le n}\|\phi(t_i)\|$.   Take
$\eta>0$;   let $\eta'>0$ be such that $(n+1)M\eta'\le\eta$ and
$\|\int_H\phi\|\le\eta$ whenever $\mu H\le(n+1)\eta'$ (see 2B(iii)).
Then we can
find a family $\langle[a_{ij},b_{ij}]\rangle_{i\le n,j\le r(i)}$ of
non-overlapping closed intervals such that
$\mu(E_i\triangle\bigcup_{j\le r(i)}[a_i,b_i])\le\eta'$ and
$t_i-\delta(t_i)\le a_{ij}\le b_{ij}\le t_i+\delta(t_i)$ for each $i\le
n$, $j\le r(i)$.   Write $t_{ij}=t_i$ for $i\le n$, $j\le r(i)$.   Then
$\langle([a_{ij},b_{ij}],t_{ij})\rangle_{i\le n,j\le r(i)}$ can be
extended to a McShane partition of $[0,1]$ subordinate to $\delta$.   So
writing $F_i=\bigcup_{j\le r(i)}[a_{ij},b_{ij}]$ for each $i$,
$F=\bigcup_{i\le n}F_i$, we have

\centerline{$\|\int_F\phi-\sum_{i\le n,j\le r(i)}(b_{ij}-
a_{ij})\phi(t_{ij})\|\le\epsilon$}

\noindent by \GoIJ, Theorem 5;  that is,

\centerline{$\|\int_F\phi-\sum_{i\le n}\mu F_i\phi(t_i)\|\le\epsilon$.}

\noindent Next,

\centerline{$\|\int_E\phi-\int_F\phi\|\le\epsilon$}

\noindent because $\mu(E\triangle F)\le(n+1)\eta$.
Also

\centerline{$\|\sum_{i\le n}\mu F_i\phi(t_i)
-\sum_{i\le n}\mu E_i\phi(t_i)\|\le M(n+1)\eta'\le \eta$.}

\noindent Putting these together,

\centerline{$\|\int_E\phi-\sum_{i\le n}\mu E_i\phi(t_i)\|\le
\epsilon+2\eta$;}

\noindent as $\eta$ is arbitrary we have the result.

\medskip

{\bf 2I Theorem} Let
$X$ be a Banach space.   Let $\langle
\phi_n\rangle_{n\in\Bbb N}$ be a sequence of McShane integrable
functions from $[0,1]$ to $X$, and suppose that
$\phi(t)=\lim_{n\to\infty}\phi_n(t)$ exists in $X$ for every $t\in
[0,1]$.
If moreover  the limit

\centerline{$\nu E =\lim_{n\to\infty}\int_E\phi_n$}

\noindent exists in $X$, for the weak topology, for every measurable
$E\subseteq[0,1]$,
$\phi$ is McShane
integrable and $\int\phi=\nu([0,1])$.

\medskip

\noindent{\bf proof}  Fix $\epsilon > 0$.   Write $\mu$ for Lebesgue
measure, $\Sigma$ for the algebra of Lebesgue measurable subsets of
$[0,1]$.

\medskip

{\bf (a)} For $t\in [0,1]$, $n\in\Bbb N$ set $q_n(t)=\sup_{j\ge i\ge
n}\|\phi_j(t)-\phi_i(t)\|$.
For each $t$,
write $r(t)=\min\{n:q_n(t)\le \epsilon,\,\|\phi(t)\|\le n\}$;  set
$A_k=\{t:r(t)=k\}$ for each $k$.   For each $k\in\Bbb N$, let
$W_k\supseteq A_k$ be a measurable set with $\mu_*(W_k\setminus A_k)=0$;
set $V_k=W_k\setminus\bigcup_{j<k}W_j$ for each $k$, so that $\langle
V_k\rangle_{k\in\Bbb N}$ is a disjoint cover of $[0,1]$ by measurable
sets,
and $A_k\subseteq\bigcup_{j\le k}V_j$ and $\mu_*(V_k\setminus A_k)=0$
for each $k$.  For each $k$, write $V^*_k=\bigcup_{j\le
k}V_j=\bigcup_{j\le k}W_j$;  take
$\eta_k>0$ such that $\|\nu E\|\le
2^{-k}\epsilon$ whenever $\mu E\le\eta_k$ (see (b) and (c) of 2F above);
let $G_k\supseteq V^*_k$
be an open set such that $\mu(G_k\setminus V^*_k)\le \min(\eta_k,
2^{-k}\epsilon)$.

\medskip

{\bf (b)} If $k\in\Bbb N$ and $E\subseteq V^*_k$ is measurable, then
$\|\nu E-\int_E\phi_k\|\le\epsilon\mu E$.   To see this, it is enough to
consider the case $E\subseteq V_j$ where $j\le k$.   In this case,
observe that

\centerline{$\|\nu E-\int_E\phi_k\|\le\limsup_{n\to\infty}
\|\int_E\phi_n-\int_E\phi_k\|\le\sup_{n\ge k}\underline{\int}_E
\|\phi_n(t)-\phi_k(t)\|\mu(dt)$}

\noindent by Lemma 2G.   Now $\mu_*(E\setminus A_j)=0$ and for $t\in
A_j$ we have
$\|\phi_n(t)-\phi_k(t)\|\le q_j(t)\le \epsilon$ for every $n\ge k$, so

\centerline{$\underline{\int}_E
\|\phi_n(t)-\phi_k(t)\|\mu(dt)\le \epsilon\mu E$}

\noindent for every $n\ge k$, giving the result.

\medskip

{\bf (c)}
For each $k\in\Bbb N$ let
$\delta_k:[0,1]\to\loibr 0,\infty[$ be a gauge such that

\centerline{$\|\int_E\phi_k
-\sum_{i\le n}\mu E_i\phi_k(t_i)\|\le 2^{-k}\epsilon$}

\noindent whenever $E_0,\ldots,E_n$ are disjoint measurable sets with
union $E$ and $t_0,\ldots,
\discrversionC{$\discretionary{}{}{}$}{}
t_n\in [0,1]$ are such that
$E_i\subseteq[t_i-\delta_k(t_i),t_i+\delta_k(t_i)]$ for each $i$;  such
a gauge exists by
Lemma
2H.   Choose $\delta:[0,1]\to\loibr 0,\infty[$ such that
$\delta(t)\le\min(\epsilon,\delta_k(t))$ and
$[0,1]\cap[t-\delta(t),t+\delta(t)]\subseteq G_k$ for $t\in A_k$.

\medskip

{\bf (d)} Let $\langle([a_i,b_i],t_i)\rangle_{i\in\Bbb N}$ be a McShane
partition of $[0,1]$ subordinate to $\delta$.   We seek to estimate
$\|\nu([0,1])-w\|$, where $w=\sum_{i\le n}(b_i-a_i)\phi(t_i)$.

Set $I_k=\{i:i\le n,\,t_i\in A_k\}$ for each $k$;  of course all but
finitely many of the $I_k$ are empty.   For $i\in I_k$, set
$E_i=[a_i,b_i]\cap V^*_k$.   We have
$[a_i,b_i]\subseteq[t_i-\delta(t_i),t_i+\delta(t_i)]\subseteq G_k$, so
$\sum_{i\in
I_k}\mu([a_i,b_i]\setminus E_i)\le 2^{-k}\epsilon$, and $\sum_{i\in
I_k}\mu([a_i,b_i]\setminus E_i)\|\phi(t_i)\|\le
2^{-k}k\epsilon$, because $\|\phi(t)\|\le k$ for $t\in A_k$.
Consequently, if we write

\centerline{$w_1=\sum_{i\le n}\mu E_i\phi(t_i)$,}

\noindent we shall have $\|w-w_1\|\le\sum_{k\in\Bbb N}
2^{-k}k\epsilon=2\epsilon$.

For each $i\le n$, let $k(i)$ be such that $t_i\in A_{k(i)}$.   Then we
have $\|\phi(t_i)-\phi_{k(i)}(t_i)\|\le \epsilon$ for each $i$.   So

\centerline{$\sum_{i\le n}\mu E_i\|\phi(t_i)-\phi_{k(i)}(t_i)\|
\le\sum_{i\le n}(b_i-a_i)\epsilon\le\epsilon$,}

\noindent because $\langle[a_i,b_i]\rangle_{i\le n}$ is
non-overlapping.   Accordingly, writing

\centerline{$w_2=\sum_{i\le n}\mu E_i\phi_{k(i)}(t_i)$,}

\noindent we have $\|w-w_2\|\le 3\epsilon$.

Set $H_k=\bigcup\{E_i:i\in I_k\}$ for each $k$.
Because $E_i\subseteq[t_i-\delta_k(t_i),t_i+\delta_k(t_i)]$ for each
$i\in I_k$, we have

\centerline{$\|\sum_{i\in I_k}\mu E_i\phi_k(t_i)-\int_{H_k}\phi_k\|
\le 2^{-k}\epsilon$.}

\noindent Consequently, writing

\centerline{$w_3=\sum_{k\in\Bbb N}\int_{H_k}\phi_k$,}

\noindent we have $\|w-w_3\|\le 5\epsilon$.

Next, for any $k$,
$H_k\subseteq V^*_k$, so
we have

\centerline{$\|\nu H_k-\int_{H_k}\phi_k\|\le\epsilon\mu H_k$,}

\noindent by (b) above.   So writing $w_4=\sum_{k\in\Bbb N}\nu H_k$ we
have $\|w_3-w_4\|\le\epsilon$ and $\|w-w_4\|\le 6\epsilon$.

If we set $H'_k=\bigcup\{[a_i,b_i]:i\in I_k\}$, then $\mu(H'_k\setminus
H_k)\le\eta_k$, so that $\|\nu H'_k-\nu H_k\|\le
2^{-k}\epsilon$, for each $k$.   Accordingly $\|w
-w_5\|\le 8\epsilon$, where

\centerline{$w_5=\sum_{k\in\Bbb N}\nu H'_k=\nu(\bigcup_{k\in\Bbb
N}H'_k)=\nu(\bigcup_{i\le n}[a_i,b_i])=\nu([0,1])$.}

\noindent As $\epsilon$ is arbitrary, $\int\phi$ exists and is equal to
$\nu([0,1])$.

\medskip

\noindent{\bf Problem} In this theorem we are supposing that $\phi(t)=\lim
_{n\to\infty}\phi_n(t)$
in the 
norm topology for every $t$.   Is it enough if $\phi(t)$ is the weak
limit of $\langle \phi_n(t)\rangle_{n\in\Bbb N}$ for every $t$?

\medskip

{\bf 2J Corollary} Let 
$X$ be a Banach space.

(a) Let $\langle\phi_n\rangle_{n\in\Bbb N}$ be a sequence of McShane
integrable functions from $[0,1]$ to $X$ such that
$\phi(t)=\lim_{n\to\infty}\phi_n(t)$ exists in $X$ for every $t\in
[0,1]$.
If

\centerline{$C=\{f\phi_n:f\in X^*,\,\|f\|\le 1,\,n\in\Bbb N\}$}

\noindent is uniformly integrable, then $\phi$ is McShane integrable.
In particular, if $\{\|\phi_n\|: n\in\Bbb N\}$ is dominated
by an integrable function, then $\phi$ is McShane integrable.

(b) Let $\phi:[0,1]\to X$ be a Pettis integrable function and $\langle
E_i\rangle_{i\in\Bbb
N}$ a cover of $[0,1]$ by measurable sets.   Suppose that
$\phi\times\chi(E_i)$ is McShane integrable for each $i$.
Then $\phi$ is McShane integrable.

\medskip

\noindent{\bf proof (a)} The point is that $\phi_n$, $\phi$ satisfy the
conditions of Theorem 2I.   To see this, take $E\in\Sigma$ and
$\epsilon>0$.  Because $C$ is uniformly integrable, there is an
$\eta>0$ such that $\int_H|g|\le\epsilon$
whenever $g\in C$ and $\mu H\le\eta$;  consequently
$\|\int_H\phi_n\|\le\epsilon$ for all $n\in\Bbb N$
whenever $H\in\Sigma$ and $\mu H\le\eta$.    Now set

\centerline{$A_n=\{t:\|\phi_i(t)
-\phi_j(t)\|\le\epsilon\enskip\forall\enskip i$, $j\ge n\}$;}

\noindent then $\langle A_n\rangle_{n\in\Bbb N}$ is an increasing
sequence with union $[0,1]$, so there is an $n$ such that $\mu^*A_n\ge 1
-\eta$.   Let $G\in\Sigma$ be such that $A_n\subseteq G$ and $\mu
G=\mu^*A_n$.   Then whenever
$i$, $j\ge n$ we have

\centerline{$\|\int_{E\cap G}\phi_i-\int_{E\cap G}\phi_j\|
\le\underline{\int}_{E\cap G}\|\phi_i(t)-\phi_j(t)\|\mu(dt)
\le\mu G\sup_{t\in A_n}\|\phi_i(t)-\phi_j(t)\|\le\epsilon$.}

\noindent Also $\|\int_{E\backslash G}\phi_i\|$ and $\|\int_{E\backslash
G}\phi_j\|$ are both less than or equal to $\epsilon$, so
$\|\int_E\phi_i-\int_E\phi_j\|\le 3\epsilon$.   This shows that
$\langle\int_E\phi_i\rangle_{i\in\Bbb N}$ is a Cauchy sequence and
therefore convergent, for every $E\in\Sigma$.   Accordingly the
conditions of 2I are satisfied and $\phi$ is McShane integrable.

\medskip

{\bf (b)} We apply 2I with $\phi_n(t)=\phi(t)$ for $t\in\bigcup_{i\le
n}E_i$, $0$ elsewhere.

\medskip

\noindent{\bf Remark} Part (a) is a version of Vitali's lemma.
Part (b) is a generalization of \GoIJ, Theorem 15.

\medskip

{\bf 2K} We now give a result connecting the McShane and
Talagrand integrals.   Recall that if $(S,\Sigma,\mu)$ is a
probability space, a set $A$ of
real-valued functions is {\bf stable} (in Talagrand's terminology) if
for
every $E\in\Sigma$, with $\mu E>0$, and all real numbers $\alpha<\beta$,
there are $m$, $n\ge 1$ such that
$\mu_{m+n}^*Z(A,E,m,n,\alpha,\beta)<(\mu E)^{m+n}$, where throughout the
rest of
paper we write $Z(A,E,I,J,\alpha,\beta)$ for

\centerline{$\{(t,u):t\in E^I,\,u\in
E^J,\,\exists\enskip f\in
A,\,f(t(i))\le\alpha\enskip\forall\enskip
i\in I,\,f(u(j))\ge\beta\enskip\forall\enskip j\in J\},$}

\noindent and $\mu^*_{m+n}$ is the ordinary product outer measure
on $S^m\times S^n$.   Now if $X$ is a Banach space, a function
$\phi:S\to X$ is {\bf properly measurable} if $\{h\phi:h\in
X^*,\,\|h\|\le
1\}$ is stable.   Talagrand proved (\TaHG, Theorem 8) that $\phi$ is
Talagrand integrable iff it is properly measurable and the upper
integral $\overline{\int}\|\phi(t)\|\mu(dt)$ is finite.

In particular, a Talagrand integrable function $\phi:S\to X$ must be scalarly
measurable for the completion of $\mu$ (\TaHD, 6-1-1).   So if $X$ is
separable, $\phi$ must be measurable for the completion of $\mu$
(\DUGG, II.1.2 or \TaHD, 3-1-3);  now as $\int\|\phi\|d\mu<\infty$,
$\phi$ is Bochner integrable.

\medskip

{\bf 2L Proposition} Let $X$ be a Banach space such that the unit ball
of $X^*$ is w$^*$-separable.
If $\phi:[0,1]\to X$ is a McShane integrable function then it is
properly measurable.

\medskip

\noindent{\bf proof} Let $w$ be the McShane integral of $\phi$.  Set
$A=\{h\phi:h\in X^*,\,\|h\|\le 1\}\subseteq\Bbb R^{[0,1]}$;  we have to
show
that $A$ is stable.    Note that because the unit ball of $X^*$ is
separable for the
w$^*$-topology on $X^*$, and the map $h\mapsto h\phi:X^*\to\Bbb
R^{[0,1]}$ is
continuous for the
w$^*$-topology on $X^*$ and the topology of pointwise convergence on
$\Bbb R^{[0,1]}$,  $A$ has a countable dense subset $A_0$.

Take a
non-negligible measurable $E\subseteq[0,1]$ and
$\alpha<\alpha'<\beta'<\beta$ in
$\Bbb R$.   For $m$, $n\ge 1$ set $H_{mn}=Z(A,E,m,n,\alpha,\beta)$,
$H'_{mn}=Z(A_0,E,m,n,\alpha',\beta')$;
then $H_{mn}\subseteq H'_{mn}$ and $H'_{mn}$ is measurable for
the usual (completed) product measure on $E^m\times E^n$.   We seek an
$m$
with $\mu_{2m}H'_{mm}<(\mu E)^{2m}$, writing $\mu$ for Lebesgue measure
on $[0,1]$ and $\mu_{2m}$ for Lebesgue measure on
$[0,1]^m\times[0,1]^m$.

Set $\epsilon={1\over
6}(\beta'-\alpha')\mu E$, and choose a function $\delta:[a,b]\to\loibr
0,\infty[$ such that $\|w-\sum_{i\le n}(b_i-a_i)\phi(t_i)\|\le\epsilon$
for every McShane partition $\langle([a_i,b_i],t_i)\rangle_{i\le n}$ of
$[a,b]$ subordinate to $\delta$.

Take $k\ge 1$ such that $\mu^*D\ge{1\over 2}\mu E$, where $D=
\{s:s\in E,\,\delta(s)\ge{1\over{k}}\}$.   Let $\langle[a_i,b_i]\rangle
_{i<m}$ be an enumeration of the intervals
of the form $[{j\over k},{{j+1}\over k}]$ which meet $D$ in a set of
positive measure;  set $G=\bigcup_{i<m}]a_i,b_i[$, so that $G$ is open
and $\mu G\ge{1\over 2}\mu E$.

Now suppose, if possible, that $\mu_{2m}H'_{mm}=(\mu E)^{2m}$.
The set

\centerline{$\{(t,u):t,\, u\in\prod_{i<m}(D\cap]a_i,b_i[)\}\subseteq
E^m\times E^m$
}

\noindent has non-zero outer measure,
so must meet $H'_{mm}$;  take $(t,u)$ in the intersection.   Write
$\Delta(s)=\{s':|s'-s|<\min(\delta(s),{1\over k})\}$ for each $s\in
[0,1]$.   Let
$\langle t(i)\rangle_{m\le i\le n}$ be a finite sequence in
$[a,b]\setminus G$ such that $[a,b]\setminus G\subseteq\bigcup_{m\le
i\le n}\Delta(t(i))$.   Because $[a,b]\setminus
G$ is itself a finite union of closed intervals, we can find a family
$\langle[a_i,b_i]\rangle_{m\le i\le n}$ of
non-overlapping closed intervals such that $[0,1]\setminus
G=\bigcup_{m\le i\le n}[a_i,b_i]$ and
$t(i)\in[a_i,b_i]\subseteq\Delta(t_i)$ for each $i$.
Now if we set $u(i)=t(i)$ for $m\le i\le n$, we see that
$\langle([a_i,b_i],t(i))\rangle_{i\le n}$ and
$\langle([a_i,b_i],u(i))\rangle_{i\le n}$ are both McShane partitions
of $[a,b]$ subordinate to $\delta$.   So we must have

\centerline{$\|\sum_{i\le n}(b_i-a_i)(\phi(t(i))-
\phi(u(i)))\|\le2\epsilon$.}

Now $(t,u)\in H'_{mm}$, so there is an $f\in A$ such that
$f(t(i))\le\alpha'$ and $f(u(i))\ge\beta'$ for every $i<m$.  $f$ is of
the form $h\phi$ for some $h$ of norm at most $1$, so

\centerline{$|\sum_{i<m}(b_i-a_i)(f(t(i))-f(u(i)))|\le 2\epsilon$.}

\noindent However, $f(t(i))\le\alpha'$ for each $i$ and
$\sum_{i<m}(b_i-a_i)
=\mu G$, so

\centerline{$\sum_{i<m}(b_i-a_i)f(t(i))\le\alpha'\mu G$;}

\noindent similarly $\sum_{i<m}(b_i-a_i)f(u(i))\ge\beta'\mu G$, and we
get

\centerline{$2\epsilon\ge(\beta'-\alpha')\mu G\ge
(\beta'-\alpha'){1\over
2}\mu E=3\epsilon$,}

\noindent which is absurd.

\medskip

{\bf 2M Corollary} Let
$X$ be a Banach space such that the unit
ball of $X^*$ is
w$^*$-separable.   If $\phi:[0,1]\to X$ is McShane integrable and
$\overline{\int}\|\phi(s)\|\mu(ds)<\infty$, then $\phi$ is Talagrand
integrable.

\medskip

\bigskip

\filename{FMp91c.tex}

\versiondate{20.1.92}

\noindent{\bf 3. Examples} 
In this section we give examples to show that the results above
are more or less complete in their own terms.   In particular, a
McShane integrable function need not be Talagrand integrable (3A, 3G)
and a Talagrand integrable function need not be McShane integrable
(3E).

\medskip

{\bf 3A Example} There is a bounded McShane integrable function
$\phi:[0,1]\to\ell^{\infty}(\frak c)$ which is not properly measurable
and therefore not Talagrand integrable.

\medskip

\noindent{\bf proof} Enumerate as $\langle H_{\xi}\rangle_{\xi<\frakc}$
the family of all Borel subsets of any power $[0,1]^n$ of the unit
interval which have positive Lebesgue measure in the appropriate
dimension.   Then
we can choose inductively a disjoint
family $\langle D_{\xi}\rangle_{\xi<\frakc}$ of finite sets such that
$D_{\xi}^{k(\xi)}\cap H_{\xi}\ne\emptyset$ for each $\xi$,
taking $k(\xi)$ such that $H_{\xi}\subseteq[0,1]^{k(\xi)}$.   Define
$\phi:[0,1]\to\ell^{\infty}(\frak c)$ by saying that $\phi(s)=e_{\xi}$
when $s\in D_{\xi}$, $0$ when
$s\in[0,1]\setminus\bigcup_{\xi<\frakc}D_{\xi}$, writing $e_{\xi}$ for
the unit vector of $\ell^{\infty}(\frak c)$ with $e_{\xi}(\eta)=1$ if
$\eta=\xi$, $0$ otherwise.

To see that $\phi$ is McShane integrable, with integral $0$,
let $\epsilon>0$.   For each $\xi<\frak c$ let $G_{\xi}\supseteq
D_{\xi}$ be a relatively open subset of $[0,1]$ of measure at most
$\epsilon$;  let $\delta(s)$ be the distance from $s$ to $[0,1]\setminus
G_{\xi}$  if $s\in D_{\xi}$, $1$ if
$s\in[0,1]\setminus\bigcup_{\xi<\frakc}D_{\xi}$.   If
$\langle([a_i,b_i],t_i)\rangle_{i\in I}$ is any McShane partition of
$[0,1]$
subordinate to $\delta$, then

\centerline{$|(\sum_{i\in J}(b_i-a_i)\phi(t_i))(\xi)|=\sum_{i\in
J,t_i\in D_{\xi}}(b_i-a_i)\le\mu G_{\xi}\le\epsilon$}

\noindent for any $J\subseteq I$, $\xi<\frak c$, so that

\centerline{$\|\sum_{i\in J}(b_i-a_i)\phi(t_i)\|\le\epsilon$}

\noindent for every finite $J\subseteq I$.   As $\epsilon$ is arbitrary,
$\phi$ is McShane integrable, with integral 0.

To see that $\phi$ is not properly measurable, set

\centerline{$A=\{h\phi:h\in(\ell^{\infty}(\frak c))^*,\,\|h\|\le 1\}$,}

\noindent and consider,
for $m$, $n\ge 1$, the set
$H_{mn}=Z(A,[0,1],m,n,0,1)$.

Suppose, if possible, that there are $m$, $n\ge 1$ such that
$(\mu )_{m+n}^*H_{mn}<1$.   In this case there is a
non-negligible measurable $H\subseteq([0,1]^m\times[0,1]^n)\setminus
H_{mn}$.   Set

\centerline{$H'=\{u:u\in[0,1]^n,\,(\mu )_n\{t:(t,u)\in H\}>0\}$;}

\noindent then $\mu_nH'>0$, so there is a $\xi<\frak c$ such that
$H_{\xi}\subseteq H'$, and a $u\in D_{\xi}^n\cap H'$.   Now
$H^{-1}[\{u\}]=\{t:(t,u)\in H\}$ is
non-negligible, and $D_{\xi}$ is finite, so there is a $t\in
H^{-1}[\{u\}]$ such that no coordinate of $t$ belongs to $D_{\xi}$.   In
this case, taking $h_{\xi}(x)=x(\xi)$ for $x\in\ell^{\infty}(\frak c)$,
and $f=h_{\xi}\phi\in A$, we see that $f(u(j))=1$ for each $j<n$, but
that $f(t(i))=0$ for each $i<m$;  so that $(t,u)\in H_{mn}$ and
$(t,u)\notin H$, which is absurd.

Thus $\phi$ is not properly measurable and therefore not Talagrand
integrable.

\medskip

\noindent{\bf Remark} Observe that $\phi$ is not measurable, and either
for this reason, or because it is not Talagrand integrable, cannot be
Bochner integrable.

\medskip

{\bf 3B} The next example will be of a Talagrand integrable function
which is not McShane integrable.   It relies on a couple of lemmas.

\medskip

\noindent{\bf Lemma} Let $(S,\Sigma,\mu)$ be a probability space and
$B\subseteq\Bbb R^S$ a stable set of functions.   Suppose that
$A\subseteq\Bbb R^S$ is a countable set of functions such that
$A\setminus U$ is stable whenever $U$ is an open subset of $\Bbb R^S$
(for the pointwise topology of $\Bbb R^S$) including $B$.   Then $A$ is
stable.

\medskip

\noindent{\bf proof} Let $E\in\Sigma$ be such that $\mu E>0$, and take
$\alpha<\beta$ in $\Bbb R$.   Let $m$, $n\ge 1$ be such
that
\linebreak ( $\mu_{m+n}^*Z(B,E,m,n,\alpha,\beta) <\mu E)^{m+n}$.
Let $G\subseteq (E^m\times E^n)\setminus Z(B,E,m,n,\alpha,\beta)$ be a
measurable set of
non-zero measure.
For any $(t,u)\in G$, we see that

\centerline{$U_{tu}=\{f:f\in\Bbb R^S,\,\exists\enskip
i<n,\,f(t(i))>\alpha\}\cup\{f:f\in\Bbb R^S,\,\exists\enskip
j<m,\,f(u(j))<\beta\}$
}

\noindent is an open set including $B$, so
$A\setminus U_{tu}$ is stable, and there are $m_{tu}\ge 1$, $n_{tu}\ge
1$ such that
\linebreak
$\mu_{m_{tu}+n_{tu}}Z(A,U_{tu},E,m_{tu},n_{tu},\alpha,\beta) <
 (\mu
E)^{m_{tu}+n_{tu}}$.   Let $p$, $q$ be so large that
$\mu_{m+n}^*G_1>0$, where

\centerline{$G_1=\{(t,u):(t,u)\in G,\,m+m_{tu}\le p,\,n+n_{tu}\le q\}$.}

\noindent
Now observe that if  $(v,w)\in Z(A,E,p,q,\alpha,\beta)$ and
$(t,u)=(v\restr m,w\restr n)\in G_1$ then we must have an $f\in A$ such
that $f(v(i))\le\alpha$ for every $i$ and $f(w(j))\ge\beta$ for every
$j$;   but in this case $f\notin U_{tu}$, so that $(v\restr p\setminus
m,w\restr q\setminus n)\in Z(A\setminus U_{tu},E,p\setminus m,q\setminus
n,\alpha,\beta)$.   Because
$p-m\ge m_{tu}$,
$q-n\ge n_{tu}$ we must have the measure of $Z(A\setminus
U_{tu},E,p\setminus m,q\setminus n,\alpha,\beta)$ strictly less than
$(\mu E)^{p-m+q-n}$.
But this means that for every $(t,u)\in G_1$,

\centerline{$\mu_{p-m+q-n}\{(v,w):(t^{\smallfrown}v,u^{\smallfrown}w)\in
Z(A,E,p,q,\alpha,\beta)\}$}

\noindent is less than $(\mu E)^{p-m+q-n}$.   Because
$Z(A,E,p,q,\alpha,\beta)$ is measurable and $\mu_{m+n}^*G_1>0$, we see
that \linebreak $\mu_{p+q}Z(A,E,p,q,\alpha,\beta)$ must be less than $(\mu
E)^{p+q}$.

As $\alpha$, $\beta$, $E$ are arbitrary, $A$ is stable.

\medskip

{\bf 3C Lemma} Let $H\subseteq\Bbb R$ be a measurable set such that
$\mu (H\cap[a,b])>0$ whenever $a<b$ in $\Bbb R$, writing $\mu $ for
Lebesgue measure.

(a) If $D$ is a
non-negligible subset of $\Bbb R$, then
$\mu ^*(H\cap(t+D))>0$ for almost all $t\in\Bbb R$.

(b) If $\langle D_n\rangle_{n\in\Bbb N}$ is a sequence of
non-negligible subsets of $\Bbb R$, then there is a sequence $\langle
t_n\rangle_{n\in\Bbb N}$ such that $t_m\in D_m$ and $t_m+t_n\in H$
whenever $m< n\in\Bbb N$.

\medskip

\noindent{\bf proof (a)} Let $E$ be a Lebesgue measurable set such that
$D\subseteq E$ and $(\mu )_*(E\setminus D)=0$.   The function

\centerline{$t\mapsto\mu ^*(H\cap(t+D))=\mu ^*((H-t)\cap D)
=\mu ((H-t)\cap E)$}

\noindent is continuous, so $F=\{t:\mu (H\cap(t+D))=0\}$ is closed.
Let $\theta$ be a
translation-invariant multiplicative lifting of Lebesgue measure
(\ITFG).   Then

\centerline{$(\theta(H)-t)\cap \theta(E)=\theta(H-t)\cap\theta(E)
=\theta((H-t)\cap E)=\emptyset$}

\noindent for every $t\in F$, so $\theta(H)\cap(\theta(E)+F)=\emptyset$.
But $\mu (\theta(H)\triangle H)=0$ so $\mu (\theta(H)\cap[a,b])>0$
whenever $a<b$ and $\theta(E)+F$ cannot include any interval.
Consequently one of $\theta(E)$, $F$ must be negligible 
(\RuGD, chap. 8, exercise 12).   But
$\mu \theta(E)=\mu E=\mu^*_LD>0$, so $\mu F=0$, as claimed.

\medskip

{\bf (b)}
Choose $t_m$, $D_{nm}$ inductively, as follows.   Start with
$D_{n0}=D_n$ for every $n$.   Given that $\mu ^*D_{nm}>0$ for all $n\ge
m$ and that $t_i+t\in H$ whenever $i<m$, $t\in\bigcup_{n\ge m}D_{nm}$,
observe that by (a) we have

\centerline{$\mu \{t:\exists\enskip n\ge m,\,\mu (D_{nm}
\cap(H-t))=0\}=0$.}

\noindent So we can find $t_m\in D_{mm}$ such that $\mu ^*(D_{nm}\cap
(H-t_m))>0$ for every $n>m$.   Set $D_{n,m+1}=D_{nm}\cap
(H-t_m)$ for $n>m$, and continue.

\medskip

{\bf 3D Lemma} Let $H\subseteq\Bbb R$ be a measurable set such that
$\mu (H\cap[a,b])>0$ whenever $a<b\in\Bbb R$.
Let $\Cal C$ be the set of
subsets $C$ of $[0,1]$ such that $s+t\notin H$ whenever $s$, $t$ are
distinct members of $C$.   Then the set
$B$ of characteristic functions of members of $\Cal C$ is stable.

\medskip

\noindent{\bf proof}
Let $E$ be a
non-negligible measurable subset of $[0,1]$ and $\alpha<\beta$.   If
$\alpha<0$ or $\beta>1$ then $Z(B,E,1,1,\alpha,\beta)=\emptyset$.   If
$0\le\alpha<\beta\le 1$ then

\centerline{$Z(B,E,1,2,\alpha,\beta)\subseteq\{(t,(u_0,u_1)):
u_0=u_1\enskip or\enskip u_0+u_1
\notin H\}$.}

\noindent But we know from Lemma 3Ca that $\mu \{u:u\in E,\,u_0+u\in
H\}>0$ for almost all $u_0$, so that
$\gamma=\mu_2\{(u_0,u_1):u_0,\,u_1\in E,\,u_0+u_1\in H\}>0$, and now
$\mu_3Z(B,E,1,2,\alpha,\beta)
\le\mu E((\mu E)^2-\gamma)<(\mu E)^3$.

\medskip

{\bf 3E Example} There is a bounded Talagrand integrable function
$\phi:
\discrversionC{$\discretionary{}{}{}$}{}
[0,1]\to\ell^{\infty}(\Bbb N)$ which is not McShane integrable.

\medskip

\noindent{\bf proof (a)} Let $H$ be an F$_{\sigma}$ subset of $\Bbb R$
such that $0<\mu (H\cap[a,b])<b-a$ whenever $a<b$ in $\Bbb R$;  let
$\langle H_n\rangle_{n\in\Bbb N}$ be an increasing sequence of closed
sets with union $H$.   Let $\Cal C$, $B\subseteq\{0,1\}^{[0,1]}$ be
defined from $H$ as in Lemma 3D, so that $B$ is stable.   For each
$n\in\Bbb N$ let $A_n$ be the countable set of
functions $f:[0,1]\to\{0,1\}$ which have total variation at most
$n$, jump only at rational points, and are such that $\min(f(s),f(t))=0$
whenever $s<t$ and $s+t\in H_n$;  set $A=\bigcup_{n\in\Bbb N}A_n$.
Then we see that in the compact space $\{0,1\}^{[0,1]}$

\centerline{$\bigcap_{n\in\Bbb N}\overline{\bigcup_{m\ge n}A_m}\subseteq
B$.}

\noindent Accordingly every neighbourhood $U$ of $B$ in $\Bbb R^{[0,1]}$
must include all but finitely many of the $A_n$.   On the other hand,
each $A_n$ is stable, being comprised of functions of variation at most
$n$.   So
Lemma 3B tells us that $A$ is stable.

Enumerate $A$ as $\langle f_n\rangle_{n\in\Bbb N}$ and define
$\phi:[0,1]\to\ell^{\infty}(\Bbb N)$ by setting $\phi(t)=\langle
f_n(t)\rangle_{n\in\Bbb N}$ for each $n\in\Bbb N$.   Then
$\{h\phi:h\in(\ell^{\infty})^*,\,\|h\|\le 1\}$ is precisely the balanced
closed convex hull of $A$, which by \TaHD, 11-1-1, is stable.   So
$\phi$ is Talagrand integrable (since of course $\|\phi(t)\|\le 1$ for
every $t$).

\medskip

{\bf (b)} But suppose, if possible, that $\phi$ is McShane integrable.
Let $w$ be its McShane integral, and let $\delta:[0,1]\to\loibr
0,\infty[$ be such that $\|w-\sum_{i\le n}(b_i-a_i)\phi(t_i)\|\le{1\over
5}$ for every McShane partition $\langle([a_i,b_i],t_i)\rangle_{i\le
n}$ of $[0,1]$ subordinate to $\delta$.   Let $k\ge 5$ be such that
$D=\{t:\delta(t)\ge{1\over k}\}$ has outer measure
at least ${4\over 5}$.   Let $\langle [a_i,b_i]\rangle_{i<m}$ enumerate
the intervals of the form $[{j\over k},{{j+1}\over k}]$ which meet $D$
in a set of positive measure;  then ${m\over k}\ge{4\over 5}$.
Set $G=\bigcup_{i<m}]a_i,b_i[$;  then (as in the proof of 2B above) we
can find
$\langle([a_i,b_i],t_i)\rangle_{m\le i\le n}$ such that $\bigcup_{m\le
i<n}[a_i,b_i]=[0,1]\setminus G$ and $t_i-\delta(t_i)\le a_i\le t_i\le
b_i\le t_i+\delta(t_i)$ for $m\le i\le n$.

Because both $H$ and its complement meet every
non-trivial interval in a set of
non-zero measure, Lemma 3Cb tells us that there are sequences
$\langle
t_i\rangle_{i<m}$, $\langle u_i\rangle_{i<m}$ such that

\centerline{$t_i\in D\cap[a_i,b_i]$, $u_i\in D\cap[a_i,b_i]$ for all
$i<m$,}

\centerline{$t_i+t_j\in H$, $u_i+u_j\notin H$ whenever $i<j<m$.}

\noindent Now if we set $u_i=t_i$ for $m\le i\le n$, both
$\langle([a_i,b_i],t_i)\rangle_{i\le n}$ and
$\langle([a_i,b_i],
\discrversionC{$\discretionary{}{}{}$}{}
u_i)\rangle_{i\le n}$ are McShane partitions of
$[0,1]$ subordinate to $\delta$.   So we must have

\centerline{${2\over 5}\ge\|\sum_{i\le n}(b_i-a_i)\phi(t_i)-\sum_{i\le
n}(b_i-a_i)\phi(u_i)\|={1\over k}\|\sum_{i<m}\phi(t_i)-\phi(u_i)\|$.}

Let $l\ge 2m$ be such that $u_i+u_j\in H_l$ whenever $i<j<m$.   Then
$t_i+t_j\notin H_l$ for
$i<j<m$, so there is an $f\in A_l$ such that $f(t_i)=
1$ for every $i<m$, while $f(u_i)=0$ for all except at most one $i<m$.
Consequently $|\sum_{i<m}f(t_i)-f(u_i)|\ge m-1$.   But $f=f_r$ for some
$r$, so $\|\sum_{i<m}\phi(t_i)-\phi(u_i)\|\ge m-1$, and
${2\over 5}\ge{{m-1}\over k}\ge{3\over 5}$, which is absurd.

Thus $\phi$ is not McShane integrable, and provides the required
example.

\medskip

\noindent{\bf Remark} Because $\phi$ is Talagrand integrable, it must be
Pettis integrable;  thus we see that the condition \lq$\phi$ is
measurable' in Gordon's theorem (1Ca) is necessary.

\medskip

{\bf 3F} We conclude with notes on two well-known examples.

\medskip

\noindent{\bf Example} There is a bounded
McShane integrable, Talagrand integrable
function
$\phi:[0,1]\to L^{\infty}([0,1])$ which is not Bochner integrable.

\medskip

\noindent{\bf proof} Let $\phi(t)$ be the equivalence class in
$L^{\infty}$ of the characteristic function of the interval $[0,t]$ for
each $t$.   As remarked in \TaHD, exercise 3-3-2, $\phi$ is Talagrand
integrable but not Bochner integrable;  as remarked in \GoIJ,
p.~567, in a slightly different context, $\phi$ is McShane integrable.

\medskip

\noindent{\bf Remarks} We observe indeed that $\phi$ above could
be called \lq Riemann integrable', as its integrability can be witnessed
by constant gauge functions.   It is easy to see that such a function
must be both McShane integrable and Talagrand integrable.

Because the unit ball of $L^{\infty}([0,1])$ is w$^*$-separable, there is
an isometric embedding of $L^{\infty}[0,1]$ in $\ell^{\infty}(\Bbb N)$
(indeed, $L^{\infty}([0,1])$ is isomorphic to $\ell^{\infty}(\Bbb N)$
-- see \LTGG, p.~111),
so there is a bounded McShane integrable, Talagrand integrable function from
$[0,1]$ to $\ell^{\infty}(\Bbb N)$ which is not Bochner integrable.

\medskip

{\bf 3G Example} There is a McShane integrable function $\phi:[0,1]\to\ell^2
(\Bbb N)$ which is not Talagrand integrable.

\medskip

\noindent{\bf proof} Let $e_n$ be the $n$th unit vector of $\ell^2(\Bbb N)$
and 
set $\phi(t)=2^n(n+1)^{-1}e_n$ for $2^{-n-1}\le t<2^{-n}$.   Then $\int\|\phi
\|=\infty$ so $\phi$ is not Talagrand integrable, but by \GoIJ, Theorem 15, 
it is
McShane integrable.

\bigskip
\filename{FMp91d.tex}

\versiondate{20.1.92}

\noindent{\bf Acknowledgement} Part of the work of this paper was done
while the first author was visiting the Universidad Nacional de
Educaci\'on a Distancia in Madrid.
The second author was partially supported by D.G.I.C.Y.T. grant
PB88-0141.
We are grateful to L.~Drewnowski for useful comments.

\bigskip

\noindent{\bf References}

\BSMHF\ J.~Bastero \& M.~San Miguel (eds.), {\it Probability and
Banach Spaces.}  Springer, 1986 (Lecture Notes in Math. 1221).

\DiHD\ J.~Diestel, {\it Sequences and Series in Banach Spaces.}
Springer, 1984 (Graduate Texts in Mathematics 94).

\DUGG\ J.~Diestel \& J.~J.~Uhl Jr., {\it Vector Measures.}   Amer.~Math.
~Soc., 1977 (Math.~Surveys 15).

\DrHF\ L.~Drewnowski, `On the Dunford and Pettis integrals', pp.~1-15
in \BSMHF.








\GoIJ\ R.~A.~Gordon, \lq The McShane integral of
Banach-valued functions', Illinois J. Math. 34 (1990)
557-567.

\ITFG\ A.~\& C.~Ionescu Tulcea, `On the existence of a lifting commuting
with the left translations of an arbitrary locally compact group',
pp.~63-97 in \LNFG.



\LNFG\ L.~M.~LeCam \& J.~Neyman (eds.) {\it Proc. Fifth Berkeley 
Symposium in Mathematical Statistics and Probability,} vol. II.
Univ. California Press, 1967.

\LTGG\ J.~Lindenstrauss \& L.~Tzafriri, {\it Classical Banach Spaces I.}
Springer, 1977.

\McSHC\ E.~J.~McShane, {\it Unified Integration.}   Academic, 1983.



\RuGD\ W.~Rudin, {\it Real and Complex Analysis.}  McGraw-Hill, 1974.

\SFpIJ\ S.~Shelah \& D.~H.~Fremlin, \lq Pointwise compact and stable sets
of measurable functions', to appear in J. Symbolic Logic.

\TaHD\ M.~Talagrand, {\it Pettis integral and measure theory.}
Mem.~Amer.~Math. Soc. 307 (1984).

\TaHG\ M.~Talagrand, \lq The
Glivenko-Cantelli problem', Annals of Probability 15 (1987)
837-870.

\bigskip

\noindent{\it Addresses:}

D.~H.~Fremlin, Mathematics Department, University of Essex, Colchester CO4
3SQ, England

J.~Mendoza, Departamento de An\'alisis Matem\'atico, Universidad
Complu\-tense de Madrid, 28040 Madrid, Spain

\end